\documentclass[twoside,8pt,B5]{article}

\usepackage{graphicx,amsmath,latexsym,amssymb,amsthm,fancyhdr}

\font\chuto=cmbx10 at 16pt \font\chudam=cmbxsl10 \font\chudams=cmbxsl8

\newtheorem{thm}{Theorem}[section]

\newtheorem{lem}[thm]{Lemma}

\theoremstyle{definition}

\numberwithin{equation}{section}
\newtheorem{condition}[thm]{Condition}

\newcommand{\R}{\ensuremath{\mathbb{R}}}
\newcommand{\E}{\ensuremath{\mathit{E}}}
\newcommand{\PP}[1]{\ensuremath{\mathit{P}}\left(#1\right)}
\newcommand{\var}[1]{\ensuremath{\mathrm{Var}}\left(#1\right)}
\newcommand{\norm}[1]{\left\Vert#1\right\Vert}
\newcommand{\abs}[1]{\left\vert#1\right\vert}
\newcommand{\set}[1]{\left\{#1\right\}}
\newcommand{\paren}[1]{\left(#1\right)}

\begin{document}

\setcounter{page}{1}

\thispagestyle{empty}

\setlength{\unitlength}{1cm}
\begin{picture}(0,0)
\put(-0.7,2.0){\chudam The 21$^{st}$ Annual Meeting in Mathematics (AMM 2016)} \put(-0.7,1.5){\chudam Annual Pure and Applied Mathematics Conference 2016 (APAM 2016)}
\put(-0.7,1.0){\chudams Department of Mathematics and Computer Science, Faculty of Science, Chulalongkorn University}
\put(-0.7,0.5){\chudams Speaker: M.~Tuntapthai}  
\put(-.7,0.2){\line(1,0){14.7}}
\put(-.7,0.17){\line(1,0){14.7}}
\put(-.7,0.15){\line(1,0){14.7}}
\end{picture}
\vspace*{-2.5cm}

\hfill \thepage

\vskip 3cm

\centerline {\bf \chuto  Bounds of the Normal Approximation for }
\medskip
\centerline {\bf \chuto  Linear Recursions with Two Effects }

\renewcommand{\thefootnote}{\fnsymbol{footnote}}

\vskip.8cm
\centerline {Mongkhon Tuntapthai{\footnote{\textit{Corresponding author}}}{\footnote{\textit{The author is supported by the young researcher development project of Khon Kaen University}}}
}

\renewcommand{\thefootnote}{\arabic{footnote}}

\vskip.5cm
\centerline{Department of Mathematics, Faculty of Science, Khon Kaen University
}
\centerline{\texttt{mongkhon@kku.ac.th}
}

\vskip .5cm

\begin{abstract}

Let $X_0$ be a non-constant random variable with finite variance. 
Given an integer $k\ge2$, define a sequence $\set{X_n}_{n=1}^\infty$ of approximately linear recursions 
with small perturbations $\set{\Delta_n}_{n=0}^\infty$ by 
\[
	X_{n+1} = \sum_{i=1}^k a_{n,i} X_{n,i} + \Delta_n 
	\quad \text{for all } n\ge0
\]
where $X_{n,1},\dots,X_{n,k}$ are independent copies of the $X_n$ and $a_{n,1},\dots,a_{n,k}$ are real numbers.
In 2004, Goldstein obtained bounds on the Wasserstein distance between the standard normal distribution and the law of $X_n$ which is in the form $C \gamma^n$ for some constants $C>0$ and  $0 < \gamma < 1$.

In this article, we extend the results to the case of two effects by studying a linear model $Z_n=X_n+Y_n$ for all $n\ge0$, 
where $\set{Y_n}_{n=1}^\infty$ is a sequence of approximately linear recursions with an initial random variable $Y_0$ and perturbations $\set{\Lambda_n}_{n=0}^\infty$, 
i.e., for some $\ell \ge2$,
\[
	Y_{n+1} = \sum_{j=1}^\ell b_{n,j} Y_{n,j} + \Lambda_n 
	\quad \text{for all } n\ge0
\]
where $Y_n$ and $Y_{n,1},\dots,Y_{n,\ell}$ are independent and identically distributed random variables 
and $b_{n,1},\dots,b_{n,\ell}$ are real numbers.
Applying the zero bias transformation in the Stein\rq s equation, we also obtain the bound for $Z_n$.  
Adding further conditions that the two models $(X_n,\Delta_n)$ and $(Y_n,\Lambda_n)$ are independent and that the difference between variance of $X_n$ and $Y_n$ is smaller than the sum of variances of their perturbation parts, our result is the same as previous work. 
\end{abstract}

\medskip




\section{Introduction and Main Theorem}


Let $Z$ be a standard normally distributed random variable and
$X_0$ a non-constant random variable with finite variance. For a positive integer $k\ge2$, we consider a sequence $\{X_n\}_{n=1}^\infty$ of approximately linear recursions with perturbations $\{\Delta_n\}_{n=0}^\infty$,
\[
	X_{n+1} = \sum_{i=1}^k a_{n,i} X_{n,i} + \Delta_n
	\quad \text{ for all } n \ge 0
\]
where the $X_n$ and $X_{n,1},\dots,X_{n,k}$ are independent and identically distributed random variables 
and $a_{n,1},\dots,a_{n,k}$ are real numbers.
For all integers $n\ge0$, we introduce some notation for the model $\paren{X_n,a_n,\Delta_n}$,
\[
	\lambda_{a,n}^2
	= \sum_{i=1}^k a_{n,i}^2,
	\;\;
	\varphi_{a,n}
	= \sum_{i=1}^k \frac{|a_{n,i}|^3}{\lambda_{a,n}^3},
	\;\;
	\var{X_n}
	= \sigma_{X,n}^2
\]
and
\[
	\widetilde{X}_n
	= \frac{X_n-\E X_n}{\sigma_{X,n}}.
\]


Arising originally from statistical physics, the approximately linear recursions are special type of hierarchical strucutres and often applied to the conductivity of random mediums.
A natural way in the classical probability theory is to study limit theorems for the distributions of these models.
A strong law of large numbers for the hierarchical structure was obtained by \cite{Wehr1997, LiRogers1999, Jordan2002}.
The central limit theorem for recursions was first introduced by \cite{WehrWoo2001} 
and the bounds to normal approximation based on the Wasserstein distance were obtained by \cite{Goldstein2004}. 
The following  two conditions were used in the last work.
\begin{condition}
\label{condition1}
For each $i=1,\dots,k$, the sequence $\{a_{n,i}\}_{n=0}^\infty$ converges to some real number $a_i$ satisfying that at least two of the $a_i$\rq s are nonzero. 
Set $\lambda_{a}^2 = \sum_{i=1}^k a_i^2$.
There exist $0 < \delta_{X,2} < \delta_{\Delta,2} < 1$ and positive constants $C_{X,2}$, $C_{\Delta,2}$ such that for all $n\ge 0$, 
\[
	\var{X_n} 
	\ge C_{X,2}^2 \lambda_{a}^{2n} \paren{1-\delta_{X,2}}^{2n},
\]
\[
	\var{\Delta_n} 
	\le C_{\Delta,2}^2 \lambda_{a}^{2n} \paren{1-\delta_{\Delta,2}}^{2n}.
\]
\end{condition}


\begin{condition}
\label{condition2}
With $\delta_{X,2}$, $\delta_{\Delta,2}$ and $\lambda_{a}$ as in the Condition \ref{condition1}, there exists $\delta_{X,4} \ge 0$ and $\delta_{\Delta,4} \ge 0$ such that
\[
	\phi_{X,\Delta,2}
	= \frac{\paren{1-\delta_{\Delta,2}} \paren{1+\delta_{X,4}}^3} {\paren{1-\delta_{X,2}}^4}  
	< 1
	\;\; \text{ and } \;\;
	\phi_{X,\Delta,4}
	= \paren{\frac{1-\delta_{\Delta,4}}{1-\delta_{X,2}} }^2 
	< 1
\]
and positive constants $C_{X,4}$, $C_{\Delta,4}$ such that
\[
	\E \paren{X_n - \E X_n}^4 
	\le C_{X,4}^4 \lambda_{a}^{4n}
	 \paren{1+\delta_{X,4}}^{4n},
\]
\[
	\E \paren{\Delta_n - \E \Delta_n}^4 
	\le C_{\Delta,4}^4 \lambda_{a}^{4n} 
	\paren{1-\delta_{\Delta,4}}^{4n} .
\]
\end{condition}

Recall that the Wasserstien distance or $L^1$-distance between two laws $F$ and $G$ is given by
\[
	\norm{F-G}_1
	= \int_{-\infty}^\infty \abs{F(t) - G(t)} \,dt.
\]
For any random variable $X$, the law or cumulative distribution function of $X$ is denoted by $\mathcal L(X)$. 

\begin{thm} \textnormal{\cite{Goldstein2004}}
\label{intro:linearrecursion:Goldstein}
Under Conditions \ref{condition1} and \ref{condition2},
there exist constants $C>0$ and $\gamma \in (0,1)$ such that
\[
	\norm{\mathcal L (\widetilde{X}_n ) - \mathcal L (Z)}_1 
	\le C \gamma^n.
\]
\end{thm}




In this article, we extend the bounds to the case of two effects.
Let $\{Z_n\}_{n=0}^\infty$ be a sequence of linear model with two effects given by
\[
	Z_n = X_n + Y_n
	\quad \text{ for all } n \ge 0
\]
where $Y_0$ is a non-degenerated random and for some integer $\ell \ge2$,
\[
	Y_{n+1} = \sum_{j=1}^\ell b_{n,j} Y_{n,j} + \Lambda_n 
	\quad \text{for all } n\ge0
\]
where $b_{n,1},\dots,b_{n,\ell}$ are real numbers, 
$Y_{n,1},\dots,Y_{n,\ell}$ are independent copy of the $Y_n$ and $\Lambda_n$ is a small perturbation.
Note that the perturbations $\Delta_n$ and $\Lambda_n$ always depend on $X_n$ and $Y_n$, respectively.
From now on, we assume that random variables from two models of recursions $(X_n,\Delta_n)$ and $(Y_n,\Lambda_n)$ are independent for all $n\ge0$, and denote 
\[
	\lambda_n^2
	= \sum_{i=1}^k a_{n,i}^2 + \sum_{j=1}^\ell b_{n,j}^2,
	\quad 
	\var{Z_n}
	= \sigma_{X,n}^2 + \sigma_{Y,n}^2
	= \sigma_{n}^2
\]
and
\[
	\widetilde Z_n
	= \frac{Z_n -\E Z_n}{\sigma_n}.
\]

The bound for linear recursions with two effects is derived by adding further assumption that the difference between variances of two models $(X_n,\Delta_n)$, $(Y_n,\Lambda_n)$, is smaller than variances of perturbations, the following is our main theorem.

\begin{thm}
\label{intro:twoeffect:mainthm}
With constants $\delta_{X,2}$, $\delta_{X,4}$, $\delta_{\Delta,2}$ 
and $\delta_{Y,2}$,  $\delta_{Y,4}$, $\delta_{\Lambda,2}$ 
as in Condition \ref{condition1} and \ref{condition2} 
for the models $\paren{X_n,\Delta_n}$ 
and $\paren{Y_n,\Lambda_n}$, suppose that 
\[
	\psi_{X,Y,\Lambda}
	= \frac{\paren{1-\delta_{\Lambda,2}} \paren{1+\delta_{X,4}}^3 }{\paren{1-\delta_{Y,2}}\paren{1-\delta_{X,2}}^3}  
	< 1
	\;\; \text{ and } \;\;
	\psi_{Y,X,\Delta}
	= \frac{\paren{1-\delta_{\Delta,2}} \paren{1+\delta_{Y,4}}^3 }{\paren{1-\delta_{X,2}}\paren{1-\delta_{Y,2}}^3}  
	< 1
\]
and that
\[
	\abs{\var{X_n} - \var{Y_n}} 
	\le \frac{\var{\Delta_n} + \var{\Lambda_n}}{\max\{\lambda_{a,n}^2,\lambda_{b,n}^2 \} } ,
\]
then there exist constants $C>0$ and $\gamma\in(0,1)$ such that
\[
	\norm{\mathcal L (\widetilde Z_n) - \mathcal L (Z)}_1 
	\le C \gamma^n.
\]
\end{thm}




\section{Auxiliary Results}
\label{auxresult}

Before proving the main theorem, we present some results for the models $(X_n,\Delta_n)$ and $(Y_n,\Lambda_n)$. 
For all $n\ge 0$, let 
\[
	r_{X,n } 
	= \frac{\lambda_n \sigma_{X,n}}{\sigma_{n+1}},
	\quad 
	r_{Y,n }
	= \frac{\lambda_n \sigma_{Y,n}}{\sigma_{n+1}} .
\] 


We begin with the bounds of $r_{X,n}$ and $r_{Y,n}$.
 
\begin{lem}
\label{auxresult:lem:boundrn}
With constants $\delta_{X,2}$, $\delta_{\Delta,2}$ and $\delta_{Y,2}$, $\delta_{\Lambda,2}$ 
as in Condition \ref{condition1} for the models 
$(X_n,\Delta_n)$ and $(Y_n,\Lambda_n)$, and suppose that 
\[
	\abs{\var{X_n} - \var{Y_n}} 
	\le \frac{\var{\Delta_n} + \var{\Lambda_n}}{\max\{\lambda_{a,n}^2,\lambda_{b,n}^2 \} } ,
\]
then for an integer $p\ge1$, there exists a positive constant $C_{r,p}$ such that
\[
	\abs{r_{X,n}^p -1}
	\le C_{r,p} \set{\paren{\frac{1-\delta_{\Delta,2}}{1-\delta_{X,2}}}^n
	+ \paren{\frac{1-\delta_{\Lambda,2}}{1-\delta_{Y,2}}}^n }
\]
and 
\[
	\abs{r_{Y,n}^p -1}
	\le C_{r,p} \set{\paren{\frac{1-\delta_{\Delta,2}}{1-\delta_{X,2}}}^n
	+ \paren{\frac{1-\delta_{\Lambda,2}}{1-\delta_{Y,2}}}^n } .
\]
\end{lem}

\proof
Following the argument of \cite[Lemma 6]{WehrWoo2001}, we consider the variances of linear model of recursions
\begin{eqnarray*}
	\sigma_{n+1}^2
	&= &
	\var{Z_{n+1}} 
	\\
	&= &
	\lambda_{a,n}^2 \var{X_n} 
	+ \lambda_{b,n}^2 \var{Y_n}
	+ \var{\Delta_n}
	+ \var{\Lambda_n}
	\\
	&= &
	\lambda_n^2 \sigma_{X,n}^2 
	+ \lambda_{b,n}^2 \set{ \var{Y_n} - \var{X_n} }
	+ \var{\Delta_n}
	+ \var{\Lambda_n} ,
\end{eqnarray*}
The triangle inequality yields
\begin{eqnarray*}
	\sigma_{n+1}
	&\le &
	\lambda_n \sigma_{X,n}
	+ \sqrt{\lambda_{b,n}^2 \abs{ \var{Y_n} - \var{X_n} }}
	+ \sqrt{\var{\Delta_n}	+ \var{\Lambda_n} } 
	\\
	&\le &
	\lambda_n \sigma_{X,n}
	+ 2 \sqrt{\var{\Delta_n}	+ \var{\Lambda_n} } .
\end{eqnarray*}
Also, we note that 
\begin{eqnarray*}
	\lambda_{a,n}^2 \sigma_{X,n}^2
	&= &
	\sigma_{n+1}^2 
	- \lambda_{b,n}^2 \set{ \var{Y_n} - \var{X_n} }
	- \var{\Delta_n}
	- \var{\Lambda_n} 
	\\
	&\le &
	\sigma_{n+1}^2 
	+ \lambda_{b,n}^2 \abs{ \var{Y_n} - \var{X_n} }
	+ \var{\Delta_n}
	+ \var{\Lambda_n} 	,
\end{eqnarray*}
which implies that 
\begin{eqnarray*}
	\lambda_n \sigma_{X,n}
	&\le &
	\sigma_{n+1}
	+ \sqrt{\lambda_{b,n}^2 \abs{ \var{Y_n} - \var{X_n} }}
	+ \sqrt{\var{\Delta_n}	+ \var{\Lambda_n} } 
	\\
	&\le &
	\lambda_n \sigma_{X,n}
	+ 2 \sqrt{\var{\Delta_n}	+ \var{\Lambda_n} } .
\end{eqnarray*}
Then there exists a constant $C_{r,1}$ such that 
\begin{eqnarray*}
	\abs{r_{X,n}-1}
	&= &
	\frac{\abs{ \lambda_n \sigma_{X,n} - \sigma_{n+1} } }{\sigma_{n+1}} 
	\\
	&\le &
	\frac{2 \sqrt{\var{\Delta_n} + \var{\Lambda_n} } }{\sigma_{n+1}} 
	\\
	&\le &
	2 \sqrt{\frac{\var{\Delta_n}}{\var{X_{n+1}}}	}
	+ 2 \sqrt{ \frac{\var{\Lambda_n}}{\var{Y_{n+1}}}	}
	\\	
	&\le &
	\frac{2C_{\Delta,2} \paren{1-\delta_{\Delta,2}}^n }
	{C_{X,2} \lambda_a \paren{1-\delta_{X,2} }^{n+1} }  
	+ \frac{2C_{\Lambda,2} \paren{1-\delta_{\Lambda,2}}^n } 
	{C_{Y,2} \lambda_b \paren{1-\delta_{Y,2} }^{n+1} }  
	\\
	&\le &
	C_{r,1} \set{\paren{\frac{1-\delta_{\Delta,2}}{1-\delta_{X,2}}}^n
	+ \paren{\frac{1-\delta_{\Lambda,2}}{1-\delta_{Y,2}}}^n } .
\end{eqnarray*}
Now, since 
\[
	\abs{r^p-1} 
	= \abs{\paren{r-1+1}^p-1} \le \sum_{j=1}^p \binom{p}{j} \abs{r-1}^j
\]
and the assumption that 
$0<\delta_{X,2}<\delta_{\Delta,2}<1$ 
and $0<\delta_{Y,2}<\delta_{\Lambda,2}<1$,
there are constants $C_{r,p}$ such that
\[
	\abs{r_{X,n}^p -1}
	\le C_{r,p} \set{\paren{\frac{1-\delta_{\Delta,2}}{1-\delta_{X,2}}}^n
	+ \paren{\frac{1-\delta_{\Lambda,2}}{1-\delta_{Y,2}}}^n }
\]
and similarly, we can see that 
\[
	\abs{r_{Y,n}^p -1}
	\le C_{r,p} \set{\paren{\frac{1-\delta_{\Delta,2}}{1-\delta_{X,2}}}^n
	+ \paren{\frac{1-\delta_{\Lambda,2}}{1-\delta_{Y,2}}}^n }
\]
for all $p=1,2,3,\dots$.
\endproof



For all $n\ge 0$, let 
\[
	U_n = U_{X,n} + U_{Y,n}
\]
where
\[
	U_{X,n+1}
	= \sum_{i=1}^k \frac{a_{n,i}}{\lambda_n} 
	\paren{ \frac{X_{n,i} - \E X_{n,i}}{\sigma_{X,n}} }
	\;\; \text{ and } \;\;
	U_{Y,n+1}
	= \sum_{j=1}^\ell \frac{b_{n,j}}{\lambda_n} 
	\paren{ \frac{Y_{n,j} - \E Y_{n,j}}{\sigma_{Y,n}} } .
\]

Next, we follow the proof of \cite[Lemma 4.1]{ChenEt2010} to prepare an inequality for the Wasserstein distance between laws of $U_n$ and its zero bias transformation.

\begin{lem}
\label{auxresult:lem:boundLU}
For all integers $n\ge1$ and the zero bias transformation
$U_n^\ast$, $\widetilde X_n^\ast$, $\widetilde Y_n^\ast$ 
of the $U_n$, $\widetilde X_n$, $\widetilde Y_n$, respectively, 
we have 
\[
	\norm{\mathcal L (U_n) - \mathcal L (U_n^\ast) }_1
	\le \norm{\mathcal L (\widetilde X_n) - \mathcal L (\widetilde X_n^\ast)}_1
	+ 	\norm{\mathcal L (\widetilde Y_n) - \mathcal L (\widetilde Y_n^\ast)}_1 .
\]
\end{lem}

\proof
Set $m=k+\ell$. Let
\[
	\xi_i =
	\begin{cases}
		\paren{ X_{n,i} - \E X_{n,i} } / \sigma_{X,n}
		& \text{ for } i = 1,\dots, k
		\\
		\paren{ Y_{n,i-k} - \E Y_{n,i-k} } / \sigma_{Y,n}
		& \text{ for } i = k+1,\dots,m
	\end{cases}
\]
and
\[
	\alpha_{n,i} =
	\begin{cases}
		a_{n,i}
		& \text{ for } i = 1,\dots, k
		\\
		b_{n,i-k}
		& \text{ for } i = k+1,\dots,m.
	\end{cases}
\]

Note that $U_{n+1}$ is a sum of independent random variables and can be written as 
\[
	U_{n+1}
	= \sum_{i=1}^m \frac{\alpha_{n,i}}{\lambda_n} \xi_i.
\]
Let $I$ be a random index independent of all other variables and satisfying that 
\[
	\PP{I = i} 
	= \frac{\alpha_{n,i}^2}{\lambda_n^2}
	\quad \text{for } i = 1,\dots,m.
\]
By the result of \cite[Lemma 2.8]{ChenEt2010}, 
the random variable
\[
	U_{n+1}^\ast
	= U_{n+1} - \frac{\alpha_{n,I}}{\lambda_n} \paren{\xi_{I}^\ast - \xi_{I} }
 \]
has the $U_{n+1}$-zero biased distribution.
By taking the dual form of the $L^1$-distance discussed in \cite{Rachev1984}, we can see that 
\[
	\norm{\mathcal L (U_{n+1}) - \mathcal L (U_{n+1}^\ast) }_1
	\; = \; \inf \E \abs{X-Y}
	\; \le \; \E \abs{U_{n+1} - U_{n+1}^\ast}
\]
where the infimum is taken over all coupling of $X$ and $Y$ having the joint distribution with $\mathcal L (U_{n+1})$ and its zero bias distribution.


Let $V_1,\dots, V_m$ be independent uniformly distributed random variables on $[0, 1]$. 
For $i=1,\dots,m$, let $\xi_i^\ast$ be the zero bias transformation of $\xi_i$. 
Let $F_\xi$ and $F_{\xi^\ast}$ be the distribution functions of $\xi$ and $\xi^\ast$, respectively.
Set
\[
	(\xi_i,\xi_i^\ast)
	= \paren{F_{\xi_i}^{-1} (V_i) , F_{\xi_i^\ast}^{-1} (V_i) }
	\quad \text{ for all }
	i=1,\dots,m.
\]
By the results of \cite{Rachev1984}, we obtain that 
\[
	\E \abs{\xi_i - \xi_i^\ast} = 
	\begin{cases}
		\norm{\mathcal L (\widetilde X_n) - \mathcal L (\widetilde X_n^\ast)}_1
		& \text{ for }  i=1,\dots,k
		\\
		\norm{\mathcal L (\widetilde Y_n) - \mathcal L (\widetilde Y_n^\ast)}_1
		& \text{ for }  i=k+1,\dots,m.
	\end{cases}
\]
Now, we obtain
\begin{eqnarray*}
	&&\hskip-.75cm
	\norm{\mathcal L (U_{n+1}) - \mathcal L (U_{n+1}^\ast) }_1
	\\
	&\le &
	\E \abs{U_{n+1} - U_{n+1}^\ast}
	\\
	&= &
	\E \sum_{i=1}^m \frac{|\alpha_{n,i}|}{\lambda_n} 
	\abs{\xi_i - \xi_i^\ast} \boldsymbol{1} \paren{I = i}
	\\
	&= &
	\sum_{i=1}^m \frac{|\alpha_{n,i}|^3}{\lambda_n^3} 
	\E  \abs{\xi_i - \xi_i^\ast} 
	\\
	&= &
	\sum_{i=1}^k \frac{|a_{n,i}|^3}{\lambda_n^3} 
	\norm{\mathcal L (\widetilde X_n) - \mathcal L (\widetilde X_n^\ast)}_1
	+ \sum_{j=1}^\ell \frac{|b_{n,j}|^3}{\lambda_n^3} 
	\norm{\mathcal L (\widetilde Y_n) - \mathcal L (\widetilde Y_n^\ast)}_1
	\\
	&= &
	\frac{\lambda_{a,n}^3 \, \varphi_{a,n}}{\lambda_n^3}
	\norm{\mathcal L (\widetilde X_n) - \mathcal L (\widetilde X_n^\ast)}_1
	+ 	\frac{\lambda_{b,n}^3 \, \varphi_{b,n}}{\lambda_n^3} 
	\norm{\mathcal L (\widetilde Y_n) - \mathcal L (\widetilde Y_n^\ast)}_1
	\\
	&\le &
	\norm{\mathcal L (\widetilde X_n) - \mathcal L (\widetilde X_n^\ast)}_1
	+ 	\norm{\mathcal L (\widetilde Y_n) - \mathcal L (\widetilde Y_n^\ast)}_1.
\end{eqnarray*}
\endproof



\section{Proof of Main Theorem}

\proof[Proof of Theorem \ref{intro:twoeffect:mainthm}]
By the results of \cite[Theorem 4.1]{ChenEt2010}, 
we can calculate the bound on $L^1$-distance by using the zero bias transformation as follows
\begin{eqnarray}
	\norm{\mathcal L (\widetilde Z_n) - \mathcal L (Z)}_1 
	& \le &
	2 \norm{\mathcal L(\widetilde Z_n) - \mathcal L(\widetilde Z_n^\ast) }_1 .
	\label{proof:eq:LZn-LZ}
\end{eqnarray}
Moreover, we can use equivalent forms of the $L^1$-distance found in \cite{Rachev1984} and given by
\[
	\norm{\mathcal L(\widetilde Z_n) - \mathcal L(\widetilde Z_n^\ast) }_1
	\; =\;  
	\sup_{h\in \mathfrak{Lip}} \abs{\E h(\widetilde Z_n) - \E h(\widetilde Z_n^\ast)}
	\; =\;  
	\sup_{f\in \mathfrak{F_{ac}}} \abs{\E f'(\widetilde Z_n) - \E f'(\widetilde Z_n^\ast)}
\]
where 
\(
	\mathfrak{Lip}
	= \set{h\colon\R\to\R : \abs{h(x)-h(y)} \le |x-y| \; \text{ for all } x,y\in \R}
\)

\noindent
and 
\(
	\mathfrak{F_{ac}}
	= \set{f\colon\R\to\R : f \text{ is absolutely continuous, } f(0)=f'(0)=0, \, f'\in \mathfrak{Lip}} .
\)



Now, we present some facts about the Stein\rq s method for normal approximation. 
For each $f \in\mathcal F$, define $h\colon\R\to\R$ by 
\[
	h(w) = f'(w) - wf(w).
\]
By the characterization of normal distribution, $\E h(Z)=0$. 
Also, we observe that 
\[
	\abs{h'(w)} 
	\; = \; 
	\abs{f''(w) - wf'(w) - f(w)}
	\; \le \;
	1 + w^2 + \frac{w^2}{2}
\]
and hence
\[
	\abs{h(w) - h(u) }
	= \abs{\int_{u}^{w} h'(t)\,dt}
	\le \abs{w-u} +\frac{1}{2}\abs{w^3-u^3}.
\]
From the definition of zero bias transformation and that $\var{\widetilde Z_{n+1}} = 1$, we have
\begin{eqnarray}
	&&\hskip-.75cm
	\abs{\E f'(\widetilde Z_{n+1}) - \E f'(\widetilde Z_{n+1}^\ast)}
	\notag\\
	&=&
	\abs{\E f'(\widetilde Z_{n+1}) - \E \widetilde Z_{n+1} f(\widetilde Z_{n+1})}
	\notag\\
	&=&
	\abs{\E h(\widetilde Z_{n+1})}
	\notag\\
	&\le &
	\abs{\E h(\widetilde Z_{n+1}) - \E h(U_{n+1}) }
	+ \abs{\E h(U_{n+1}) }
	\notag\\
	&\le &
	\E \abs{\widetilde Z_{n+1} - U_{n+1}} 
	+ \frac{1}{2} \E \abs{\widetilde Z_{n+1}^3 - U_{n+1}^3}
	+ \abs{\E h(U_{n+1}) }
	\notag\\
	&= &
	\beta_n
	+ \abs{\E f'(U_{n+1}) - \E f'(U_{n+1}^\ast)}
	\notag\\
	&\le &
	\beta_n
	+ \norm{\mathcal L (U_{n+1}) - \mathcal L (U_{n+1}^\ast) }_1
	\notag\\
	&\le &
	\beta_n
	+ \norm{\mathcal L (\widetilde X_{n}) - \mathcal L (\widetilde X_{n}^\ast) }_1
	+ \norm{\mathcal L (\widetilde Y_{n}) - \mathcal L (\widetilde Y_{n}^\ast) }_1
	\label{proof:eq:EfZ-EfZstar}
\end{eqnarray}
where we apply Lemma \ref{auxresult:lem:boundLU} in the last inequality and denote for all $n\ge0$, 
\begin{eqnarray}
	\beta_n
	= \E \abs{\widetilde Z_{n+1} - U_{n+1}} 
	+ \frac{1}{2} \E \abs{\widetilde Z_{n+1}^3 - U_{n+1}^3} .
	\label{proof:eq:betan}
\end{eqnarray}
By (\ref{proof:eq:LZn-LZ}) and taking the supremum of (\ref{proof:eq:EfZ-EfZstar}) over $f \in \mathfrak{F_{ac}}$, we obtain
\begin{eqnarray*}
	\norm{\mathcal L (\widetilde Z_{n+1}) - \mathcal L (Z) }_1
	&\le &
	2 \norm{\mathcal L (\widetilde Z_{n+1}) - \mathcal L (\widetilde Z_{n+1}^\ast)}_1
	\\
	&\le &
	2 \beta_n
	+ 2 \norm{\mathcal L (\widetilde X_{n}) - \mathcal L (\widetilde X_{n}^\ast) }_1
	+ 2 \norm{\mathcal L (\widetilde Y_{n}) - \mathcal L (\widetilde Y_{n}^\ast) }_1 .
\end{eqnarray*}
Applying the Condition \ref{condition1} and \ref{condition2} 
for the models $(X_n,\Delta_n)$ and $(Y_n,\Lambda_n)$ in Theorem \ref{intro:linearrecursion:Goldstein}, 
there exist  positive constants $C_{X,a,\Delta}$, $C_{Y,b\Lambda}$ 
and $\gamma_{X,a,\Delta}\in(0,1)$, $\gamma_{Y,b,\Lambda}\in(0,1)$ such that for all $n\ge 0$,
\[
	\norm{\mathcal L (\widetilde X_{n}) - \mathcal L (\widetilde X_{n}^\ast) }_1
	\le C_{X,a,\Delta} \paren{\gamma_{X,a,\Delta} }^n
\]
and
\[
	\norm{\mathcal L (\widetilde Y_{n}) - \mathcal L (\widetilde Y_{n}^\ast) }_1
	\le C_{Y,b,\Lambda} \paren{\gamma_{Y,b,\Lambda} }^n .
\]
We remain to show that $\beta_n \le C_\beta \gamma_\beta^n$ for some $C_\beta>0$ and $\gamma_\beta\in(0,1)$ and the proof is completed by choosing $C=C_{X,a,\Delta} + C_{Y,b,\Lambda} + C_\beta$ and $\gamma=\max\set{\gamma_{X,a,\Delta},\gamma_{Y,b,\Lambda},\gamma_\beta}$.



\medskip

Recalling the definition of  $r_{X,n}$, $r_{Y,n}$ and $U_{X,n}$, $U_{Y,n}$ in Lemma  \ref{auxresult:lem:boundrn} and \ref{auxresult:lem:boundLU}, respectively, the linear model of recursions can be written as
\begin{eqnarray*}
	\widetilde Z_{n+1}
	&=& 
	\frac{Z_{n+1} - \E Z_{n+1}}{\sigma_{n+1}}
	\notag\\
	&=& 
	\frac{X_{n+1} - \E X_{n+1}}{\sigma_{n+1}}
	+ \frac{Y_{n+1} - \E Y_{n+1}}{\sigma_{n+1}}
	\notag\\
	&=&
	\frac{\sigma_{X,n}}{\sigma_{n+1}} \set{ \sum_{i=1}^k a_{n,i} \paren{\frac{X_{n,i} - \E X_{n,i} }{\sigma_{X,n} } }  
	+ \frac{\Delta_n - \E \Delta_n}{\sigma_{X,n}} }
	\notag\\
	&&
	+ \frac{\sigma_{Y,n}}{\sigma_{n+1}} \set{ \sum_{j=1}^\ell b_{n,j} \paren{\frac{Y_{n,j} - \E Y_{n,j} }{\sigma_{Y,n}} } 
	+ \frac{\Lambda_n - \E \Lambda_n}{\sigma_{Y,n}} }
	\notag\\
	&=&
	r_{X,n} {U}_{X,n+1} 
	+ r_{Y,n} {U}_{Y,n+1} 
	+ \Gamma_n
\end{eqnarray*}
where 
$\Gamma_n = \Gamma_{X,\Delta,n} + \Gamma_{Y,\Lambda,n}$, 
\[ 
	\Gamma_{X,\Delta,n}
	= \frac{\sigma_{X,n}}{\sigma_{n+1}} \paren{\frac{\Delta_n - \E \Delta_n}{\sigma_{X,n}} }
	\;\; \text{and } \;\;
	\Gamma_{Y,\Lambda,n}
	= \frac{\sigma_{Y,n}}{\sigma_{n+1}} \paren{\frac{\Delta_n - \E \Delta_n}{\sigma_{Y,n}} } .
\]
Using Conditions \ref{condition1} and \ref{condition2} for the models $(X_n,\Delta_n)$ and $(Y_n,\Lambda_n)$, the result of \cite[Lemma 6]{WehrWoo2001} gives that the limits 
\[
	\lim_{n\to \infty} \frac{\sigma_{X,n}}{\lambda_{a,0} \dots \lambda_{a,n-1}} 
	\;\; \text{and} \;\;
	\lim_{n\to \infty} \frac{\sigma_{Y,n}}{\lambda_{b,0} \dots \lambda_{b,n-1}} 
\]
exist in $(0,1)$,
so we have 
\[
	\lim_{n\to\infty} \frac{\sigma_{X,n+1}}{\sigma_{X,n}} = \lambda_a
	\;\; \text{ and } \;\;
	\lim_{n\to\infty} \frac{\sigma_{Y,n+1}}{\sigma_{Y,n}} = \lambda_b .
\]
Therefore, there exist  positive constants $C_{\Gamma,X,\Delta,2}$ and $C_{\Gamma,Y,\Lambda,2}$ such that 
\[
	\E \Gamma_{X,\Delta,n}^2
	\le \paren{\frac{\sigma_{X,n}}{\sigma_{X,n+1}}}^2 \frac{\var{\Delta_n}}{\var{X_n}} 
	\le C_{\Gamma,X,\Delta,2}^2 \paren{\frac{1-\delta_{\Delta,2}}{1-\delta_{X,2}}}^{2n}
\]
\[
	\E \Gamma_{Y,\Lambda,n}^2
	\le \paren{\frac{\sigma_{Y,n}}{\sigma_{Y,n+1}}}^2 \frac{\var{\Lambda_n}}{\var{Y_n}} 
	\le C_{\Gamma,Y,\Lambda,2}^2 \paren{\frac{1-\delta_{\Lambda,2}}{1-\delta_{Y,2}}}^{2n}.
\]
Moreover, there exist positive constants $C_{\Gamma,X,\Delta,4}$ and $C_{\Gamma,Y,\Lambda,4}$ such that 
\[
	\E \Gamma_{X,\Delta,n}^4
	\le \paren{\frac{\sigma_{X,n}}{\sigma_{X,n+1}}}^4 
	\E \paren{\frac{\Delta_n - \E \Delta_n}{\sigma_{X,n}} }^4
	\le C_{\Gamma,X,\Delta,4}^4 \paren{\frac{1-\delta_{\Delta,4}}{1-\delta_{X,2}}}^{4n}
\]
\[
	\E \Gamma_{Y,\Lambda,n}^4
	\le \paren{\frac{\sigma_{Y,n}}{\sigma_{Y,n+1}}}^4 
	\E \paren{\frac{\Lambda_n - \E \Lambda_n}{\sigma_{Y,n}} }^4
	\le C_{\Gamma,Y,\Lambda,4}^4 \paren{\frac{1-\delta_{\Lambda,4}}{1-\delta_{Y,2}}}^{4n} .
\]
%
%
By independence for $X_{n,i}$\rq s and $Y_{n,j}$\rq s, there exist positive constants $C_{U,X}$ and $C_{U,Y}$ such that 
\begin{eqnarray*}
	\E U_{X,n+1}^2
	&= &
	\frac{\lambda_{a,n}^2}{\lambda_n^2}
	\E \paren{\frac{X_n - \E X_n}{\sigma_{X,n}} }^2
	\le 1
	\\
	\E U_{Y,n+1}^2
	&= &
	\frac{\lambda_{b,n}^2}{\lambda_n^2}
	\E \paren{\frac{Y_n - \E Y_n}{\sigma_{Y,n}} }^2
	\le 1
	\\
	\E U_{X,n+1}^4
	&\le &
	8 \sum_{i=1}^k \frac{a_{n,i}^4}{\lambda_{n}^4} 
	\E \paren{\frac{X_n - \E X_n}{\sigma_{X,n}} }^4
	\le C_{U,X}^4 \paren{\frac{1+\delta_{X,4}}{1-\delta_{X,2}}}^{4n}
	\\
	\E U_{Y,n+1}^4
	&\le &
	8 \sum_{j=1}^\ell \frac{b_{n,j}^4}{\lambda_{n}^4} 
	\E \paren{\frac{Y_n - \E Y_n}{\sigma_{Y,n}} }^4
	\le C_{U,Y}^4 \paren{\frac{1+\delta_{Y,4}}{1-\delta_{Y,2}}}^{4n} .
\end{eqnarray*}
From Lemma \ref{auxresult:lem:boundrn} and Condition \ref{condition1} and \ref{condition2}, 
the following results will be often used for all $n\ge0$ and $p=1,2,3$,
\begin{eqnarray}
	\abs{r_{X,n}^p -1}
	\le C_{r,p} \paren{\phi_{X,\Delta,2}^n + \phi_{Y,\Lambda,2}^n}
	\label{proof:eq:rXnp-1} 
	\\
	\abs{r_{Y,n}^p -1}
	\le C_{r,p} \paren{\phi_{X,\Delta,2}^n + \phi_{Y,\Lambda,2}^n} .
	\label{proof:eq:rYnp-1}
\end{eqnarray}



Now, considering the first term of $\beta_n$ in (\ref{proof:eq:betan}), 
\begin{eqnarray*}
	&&\hskip-.75cm
	\E \abs{\widetilde Z_{n+1} - U_{n+1}} 
	\\
	&=& 
	\E \abs{\paren{r_{X,n}-1} U_{X,n+1}  
	+ \paren{r_{Y,n} - 1} U_{Y,n+1}  
	+ \Gamma_{X,\Delta,n} 
	+ \Gamma_{Y,\Lambda,n} }
	\\
	&\le &
	\abs{r_{X,n} - 1} \sqrt{ \E U_{X,n+1}^2 } 
	+ \abs{r_{Y,n} - 1} \sqrt{ \E U_{Y,n+1}^2 } 
	+ \sqrt{\E \Gamma_{X,\Delta,n}^2}
	+ \sqrt{\E \Gamma_{Y,\Lambda,n}^2}
	\\
	&\le &
	2 C_{r,1} \paren{\phi_{X,\Delta,2}^n + \phi_{Y,\Lambda,2}^n} 
	+ C_{\Gamma,X,\Delta,2} \paren{\frac{1-\delta_{\Delta,2}}{1-\delta_{X,2}} }^n
	+ C_{\Gamma,Y,\Lambda,2} \paren{\frac{1-\delta_{\Lambda,2}}{1-\delta_{Y,2}} }^n
	\\
	&\le &
	C_0 \paren{\phi_{X,\Delta,2}^n + \phi_{Y,\Lambda,2}^n} .
\end{eqnarray*}
%
%
%
%
%
%
%
%

\noindent
For the second term of $\beta_n$, 
\begin{eqnarray*}
	&& \hskip-.75cm
	\E \abs{\widetilde Z_{n+1}^3 - U_{n+1}^3}
	\\
	&= &
	\E \abs{ \paren{ r_{X,n} U_{X,n+1} 
	+ r_{Y,n} U_{Y,n+1} 
	+ \Gamma_n }^3 
	- U_{n+1}^3 }
	\\
	&= &
	\E \left| \paren{ r_{X,n} U_{X,n+1} 
	+ r_{Y,n} U_{Y,n+1} }^3 
	+ 3 \paren{ r_{X,n} U_{X,n+1} + r_{Y,n} U_{Y,n+1} }^2 \Gamma_n  \right. 
	\\
	&& \quad
	\left. \phantom{\bigg|}
	+ 3 \paren{ r_{X,n} U_{X,n+1} 
	+ r_{Y,n} U_{Y,n+1} } \Gamma_n^2 
	+  \Gamma_n^3
	- U_{n+1}^3 \right|
\\
	&\le &
	\E \abs{ \paren{ r_{X,n} U_{X,n+1} 
	+ r_{Y,n} U_{Y,n+1} }^3 - U_{n+1}^3 }
	+ 3 \E \abs{\paren{ r_{X,n} U_{X,n+1} 
	+ r_{Y,n} U_{Y,n+1} }^2 \Gamma_n }  
	\\
	&&
	+ 3 \E \abs{\paren{ r_{X,n} U_{X,n+1} 
	+ r_{Y,n} U_{Y,n+1} } \Gamma_n^2 } 
	+ \E \abs{\Gamma_n}^3
	\\
	& := &
	A_1 + A_2 + A_3 +A_4 .
\end{eqnarray*}
%

Notice that 
\begin{eqnarray*}
	A_1 
	&= &
	\E \abs{ \paren{ r_{X,n} U_{X,n+1} 
	+ r_{Y,n} U_{Y,n+1} }^3 - \paren{U_{X,n+1}+U_{Y,n+1}}^3 }
	\\
	&= &
	\E \left| \paren{r_{X,n}^3 - 1} U_{X,n+1}^3 
	+ 3\paren{r_{X,n}^2r_{Y,n}-1} U_{X,n+1}^2 U_{Y,n+1} 
	\right.
	\\
	&& \quad \;\;
	\left.
	+ 3\paren{r_{X,n}r_{Y,n}^2-1} U_{X,n+1} U_{Y,n+1}^2 
	+ \paren{r_{Y,n}^3 - 1} U_{Y,n+1}^3  \right|
	\\
	&\le &
	\E \abs{\paren{r_{X,n}^3 - 1} U_{X,n+1}^3 }
	+ 3 \E \abs{\paren{r_{X,n}^2r_{Y,n} - r_{Y,n} + r_{Y,n} -1}  U_{X,n+1}^2 U_{Y,n+1} } 
	\\
	&& 
	+ 3 \E \abs{\paren{r_{X,n} r_{Y,n}^2 - r_{X,n} + r_{X,n} -1} U_{X,n+1} U_{Y,n+1}^2 }
	+ \E \abs{\paren{r_{Y,n}^3 - 1} U_{Y,n+1}^3}
	\\
	&\le &
	\E \abs{\paren{r_{X,n}^3 - 1} U_{X,n+1}^3 }
	+ \E \abs{\paren{r_{Y,n}^3 - 1} U_{Y,n+1}^3}
	\\
	&&
	+ 3 \E \abs{\paren{r_{X,n}^2-1} r_{Y,n} U_{X,n+1}^2 U_{Y,n+1} }
	+ 3 \E \abs{\paren{r_{Y,n}-1} U_{X,n+1}^2 U_{Y,n+1} }
	\\
	&& 
	+ 3 \E \abs{\paren{r_{Y,n}^2-1} r_{X,n} U_{X,n+1} U_{Y,n+1}^2 }
	+ 3 \E \abs{\paren{r_{X,n}-1} U_{X,n+1} U_{Y,n+1}^2 }
	\\
	&\le &
	\abs{r_{X,n}^3-1} \paren{\E U_{X,n}^4}^{3/4}
	+ \abs{r_{Y,n}^3-1} \paren{\E U_{Y,n}^4}^{3/4}
	\\
	&&
	+ 3 \abs{r_{X,n}^2-1} r_{Y,n} \E U_{X,n}^2 \sqrt{\E U_{Y,n}^2}
	+ 3 r_{X,n} \abs{r_{Y,n}^2-1} \sqrt{\E U_{X,n}^2} \E U_{Y,n}^2 
	\\
	&&
	+ 3 \abs{r_{Y,n}-1} \E U_{X,n}^2 \sqrt{\E U_{Y,n}^2}
	+ 3 \abs{r_{X,n}-1} \sqrt{\E U_{X,n}^2} \E U_{Y,n}^2 
	\\
	&\le &
	8^{3/4} C_{r,3} C_{X,4}^3 \set{ \paren{\frac{\paren{1-\delta_{\Delta,2}}\paren{1+\delta_{X,4}}^3}{\paren{1-\delta_{X,2}}^4} }^n	
	+ \paren{\frac{\paren{1-\delta_{\Lambda,2}}\paren{1+\delta_{X,4}}^3}{\paren{1-\delta_{Y,2}}\paren{1-\delta_{X,2}}^3} }^n	}
	\\
	&&
	+ 8^{3/4} C_{r,3} C_{Y,4}^3 \set{ \paren{\frac{\paren{1-\delta_{\Delta,2}}\paren{1+\delta_{Y,4}}^3}{\paren{1-\delta_{X,2}}\paren{1-\delta_{Y,2}}^3} }^n	
	+ \paren{\frac{\paren{1-\delta_{\Lambda,2}}\paren{1+\delta_{Y,4}}^3}{\paren{1-\delta_{Y,2}}^4} }^n }	
	\\
	&&
	+ 6 C_{r,2} \paren{1+2C_{r,1}} 
	\paren{\phi_{X,\Delta,2}^n + \phi_{Y,\Lambda,2}^n} 
	+ 6 C_{r,1} \paren{\phi_{X,\Delta,2}^n + \phi_{Y,\Lambda,2}^n} 
	\\
	&\le &
	C_1 \paren{\phi_{X,\Delta,2}^n + \phi_{Y,\Lambda,2}^n 
	+ \psi_{X,Y,\Lambda}^n + \psi_{Y,X,\Delta}^n }.
\end{eqnarray*}

As a special case of (\ref{proof:eq:rXnp-1}) and (\ref{proof:eq:rYnp-1}) when $p=1$, we can see that for all $n\ge0$
\[
	r_{X,n} 
	\le 1 + C_{r,1} \paren{\phi_{X,\Delta,2}^n + \phi_{Y,\Lambda,2}^n}
	\le 1 + 2 C_{r,1} 
\]
\[
	r_{Y,n} 
	\le 1 + C_{r,1} \paren{\phi_{X,\Delta,2}^n + \phi_{Y,\Lambda,2}^n}
	\le 1 + 2 C_{r,1}.
\]
So, we have that 
\begin{eqnarray*}
	A_2 
	&= &
	3 \E \abs{\paren{ r_{X,n} U_{X,n+1} 
	+ r_{Y,n} U_{Y,n+1} }^2 \paren{\Gamma_{X,\Delta,n}+\Gamma_{Y,\Lambda,n} } }  
	\\
	&\le &
	6 \E \abs{\paren{ r_{X,n}^2 U_{X,n+1}^2
	+ r_{Y,n}^2 U_{Y,n+1}^2 } \paren{\Gamma_{X,\Delta,n}+\Gamma_{Y,\Lambda,n} } } 
	\\
	&\le &
	6 r_{X,n}^2 \sqrt{\E U_{X,n+1}^4 \E \Gamma_{X,\Delta,n}^2 }
	+ 6 r_{X,n}^2 \E U_{X,n+1}^2 \sqrt{\E \Gamma_{Y,\Lambda,n}^2} 
	\\
	&&
	+ 6 r_{Y,n}^2 \E U_{Y,n+1}^2 \sqrt{\E \Gamma_{X,\Delta,n}^2}
	+ 6 r_{Y,n}^2 \sqrt{\E U_{Y,n+1}^4 \E \Gamma_{Y,\Lambda,n}^2 } 
	\\
	&\le &
	12\sqrt 2 \paren{1+C_{r,1}}^2 C_{U,X}^2 C_{\Gamma,X,\Delta,2} \paren{\frac{\paren{1-\delta_{\Delta,2}}\paren{1+\delta_{X,4}}^3}{\paren{1-\delta_{X,2}}^4} }^n 
	\\
	&&
	+ 6 \paren{1+C_{r,1}}^2 C_{\Gamma,X,\Delta,2} \paren{\frac{1-\delta_{\Delta,2}}{1-\delta_{X,2}} }^n 
	+ 6 \paren{1+C_{r,1}}^2 C_{\Gamma,Y,\Lambda,2} \paren{\frac{1-\delta_{\Lambda,2}}{1-\delta_{Y,2}} }^n 
	\\
	&&
	+ 12\sqrt 2 \paren{1+C_{r,1}}^2 C_{U,Y}^2 C_{\Gamma,Y,\Lambda,2} \paren{\frac{\paren{1-\delta_{\Lambda,2}}\paren{1+\delta_{Y,4}}^3}{\paren{1-\delta_{Y,2}}^4} }^n 	
	\\
	&\le &
	C_2 \paren{\phi_{X,\Delta,2}^n + \phi_{Y,\Lambda,2}^n} 
\end{eqnarray*}
and that 
\begin{eqnarray*}
	A_3
	&= &
	3 \E \abs{\paren{ r_{X,n} U_{X,n+1} 
	+ r_{Y,n} U_{Y,n+1} } \paren{\Gamma_{X,\Delta,n}+\Gamma_{Y,\Lambda,n}}^2 }  
	\\
	&\le &
	6 \E \abs{\paren{ r_{X,n} U_{X,n+1} + r_{Y,n} U_{Y,n+1} } 
	\paren{\Gamma_{X,\Delta,n}^2+\Gamma_{Y,\Lambda,n}^2 } }  
	\\
	&\le &
	6 r_{X,n} \sqrt{ \E U_{X,n+1}^2 \E \Gamma_{X,\Delta,n}^4 }  
	+ 6 r_{X,n} \sqrt{\E U_{X,n+1}^2} \E \Gamma_{Y,\Lambda,n}^2 
	\\
	&&
	+ 	6 r_{Y,n} \sqrt{\E U_{Y,n+1}^2} \E \Gamma_{X,\Delta,n}^2
	+ 6 r_{Y,n} \sqrt{ \E U_{Y,n+1}^2 \E \Gamma_{Y,\Lambda,n}^4}  
	\\
	&\le &
	6 \paren{1+C_{r,1}} \paren{ C_{\Gamma,X,\Delta,4}^2 + C_{\Gamma,X,\Delta,2}^2 }
	\paren{\frac{1-\delta_{\Delta,2}}{1-\delta_{X,2}} }^{2n}
	\\
	&&
	+ 6 \paren{1+C_{r,1}} \paren{ C_{\Gamma,Y,\Lambda,2}^2 + C_{\Gamma,Y,\Lambda,4}^2 }
	\paren{\frac{1-\delta_{\Lambda,2}}{1-\delta_{Y,2}} }^{2n}
	\\
	&\le &
	C_3 \paren{\phi_{X,\Delta,2}^{2n} + \phi_{Y,\Lambda,2}^{2n}} .
\end{eqnarray*}
Lastly, 
\begin{eqnarray*}
	A_4
	&\le &
	\paren{\E \Gamma_n^4}^{3/4}
	\\
	&\le &
	8^{3/4} \paren{\E \Gamma_{X,\Delta,n}^4 
	+ \E \Gamma_{Y,\Lambda,n}^4 }^{3/4}	
	\\
	&\le &
	8^{3/4}  \set{ C_{\Gamma,X,\Delta,4}^4
	\paren{ \frac{1-\delta_{\Delta,4}}{1-\delta_{X,2}} }^{4n}
	+ C_{\Gamma,Y,\Lambda,4}^4
	\paren{ \frac{1-\delta_{\Lambda,4}}{1-\delta_{Y,2}} }^{4n} }^{3/4}
	\\
	&\le &
	C_4 \paren{ \phi_{X,\Delta,4}^{3n/2} 
	+ \phi_{Y,\Lambda,4}^{3n/2} }.
\end{eqnarray*}
Setting $\gamma_\beta = \max\set{\phi_{X,\Delta,2},\phi_{Y,\Lambda,2},\phi_{X,\Delta,4}^{3/2},\phi_{Y,\Lambda,4}^{3/2},\psi_{X,Y,\Lambda},\phi_{Y,X,\Delta} } \in(0,1)$

\noindent
and $C_\beta = 2C_0+4C_1+2C_2+2C_3+2C_4$, 
we obtain the claim for $\beta_n$.
\endproof



\newpage




\end{document}